\documentclass{amsart}
\usepackage{amssymb}   
\usepackage{amsmath}
\usepackage{amsthm}
\usepackage{xypic}

\def\mod{{\rm mod}\,}
\newcommand{\Pic}{\operatorname{Pic}}
\newcommand{\Bl}{\operatorname{Bl}}

\renewcommand{\O}{{\mathcal O}}

\newcommand{\Proj}{{\mathbb P}}

\newcommand{\Z}{{\mathbb Z}}

\newcommand{\C}{{\mathbb C}}
\newcommand{\p}{{\mathbb P}}

\newcommand{\map}{\dasharrow}

\newcommand{\defect}{\operatorname{def}}

\def\leq{\leqslant}
\def\geq{\geqslant}
\def\deft{~}

\def\phi{\varphi}

\theoremstyle{plain}
\newtheorem{Theorem}{Theorem}[section]
\newtheorem{Lemma}[Theorem]{Lemma}
\newtheorem{Proposition}[Theorem]{Proposition}

 \theoremstyle{definition}
\newtheorem{Definition}[Theorem]{Definition}
\newtheorem{Remark}[Theorem]{Remark}
\newtheorem{Examples}[Theorem]{Examples}

\begin{document}

\title[Varieties with quadratic entry locus, {\it II}]{Varieties with quadratic entry locus, {\it II}}
\author[Paltin Ionescu]{Paltin Ionescu}
\address{University of Bucharest and  Institute of Mathematics of the Romanian Academy}
\curraddr{Dipartimento di Matematica\\
Universit\` a degli Studi di Genova\\
Via Dodecaneso, 35\\
16146 Genova\\ Italy}
\email{ionescu@dima.unige.it}

\author[Francesco Russo]{Francesco Russo}
\address{Dipartamento di Matematica e Informatica\\
Universit\'a di Catania\\
Viale A. Doria 6\\
95125 Catania\\ Italy}
\email{frusso@dmi.unict.it}
\thanks{The first author is partially  supported  by the Italian Programme ``Incentivazione alla mobilit\` a di studiosi stranieri e italiani residenti all'estero". The second author is partially  supported  by CNPq (Centro Nacional de Pesquisa), Grants 300961/2003-0, 308745/2006-0 and 474475/2006-9, and by PRONEX/FAPERJ--Algebra Comutativa e Geometria Algebrica.}

\begin{abstract}
We continue the study, begun in \cite{QEL}, of secant defective manifolds having ``simple entry loci". We prove that such manifolds are rational and describe them in terms of tangential projections. Using also \cite{IR}, their classification  is reduced to the case of Fano manifolds of high index, whose Picard group is generated by the hyperplane section class. Conjecturally, the former should be linear sections of rational homogeneous manifolds. We also provide evidence that the classification of linearly normal dual defective manifolds with Picard group generated by the hyperplane section should follow along the same lines.  
\end{abstract}

\maketitle

\section*{Introduction}

An $n$-dimensional closed submanifold $X\subset \p^N$ is {\it secant defective} if its secant variety $SX\subseteq \p^N$ has dimension less than $2n+1$, the expected one. 
The secant defect of $X$ is then $\delta=\delta (X) := 2n+1-\dim (SX)$. Secant defective manifolds naturally fall into two categories.  
Manifolds admitting ``non-trivial projections" are those for which $SX \not = \p^N$, see \cite{Sev, Zak}, while manifolds ``of small codimension" correspond to the other case, $SX=\p^N$, see \cite{BL}. 
Let $x,y \in X$ be two general points and let $p$ be a general point on the line $\langle x,y \rangle$. Consider the closure of the locus of secants to $X$ passing through $p$. Its trace on $X$, denoted $\Sigma_p$, is called the {\it entry locus} (with respect to $p$) and has dimension $\delta$. 
Our aim is to get classification results for secant defective manifolds whose entry loci are simple enough. Consider embedded manifolds $X \subset \p^N$ as above, such that through two general points $x,y \in X$ there passes an $r$-dimensional quadric hypersurface, say $Q^r$, contained in $X$. Observe that $Q^r \subseteq \Sigma_p$ if $p \in \langle x,y \rangle$; in particular $r \leq \delta$, see \cite{KS}. When $r=1$,
we call such manifolds {\it conic-connected} (CCM for short). The extremal case $r=\delta$ was called ``manifolds with {\it local quadratic entry locus}" (abbreviated as LQELM), while the special case when $Q^r=\Sigma_p$ was named ``manifolds with {\it quadratic entry locus}" (QELM); see \cite{KS, QEL, IR}.
Being rationally connected, these special classes of secant defective manifolds may be studied in the context of Mori Theory, see \cite{Mori, Mori1, KMM} and also \cite{Debarre, Kollar, Hwang}. The interested reader can find further motivation and various examples in the introduction to \cite{QEL}. The present paper continues the line of investigation started in \cite{QEL}, see also \cite{IR}. We acknowledge once more our intellectual debt to the classical work by G. Scorza, see \cite{Scorza1, Scorza2}. The content of the paper is described below.

Tangential projections turned out to be a basic tool in the classical works (by Severi, Scorza, Terracini, and others) on secant defective manifolds. 
In what follows, by the $\delta$-partial tangential projection we mean projection from a $\delta$-codimen\-sional linear subspace (passing through $x$) of the projective tangent space at $x\in X$.  
In Theorem\deft\ref{rationality} we show that the generic $\delta$-partial tangential projection of a QELM is birational. The proof is based on a degeneration technique introduced in \cite{CMR} (see also \cite{CR}). 
The key point is that the general fibre of the (full) tangential projection is a degeneration of the general entry locus.
Conversely, if the generic $\delta$-partial tangential projection of $X$ is birational, then $X$ is an LQELM.  
Moreover, we show that defective LQELM are rational, extending the theorem on the rationality of ``manifolds with one apparent double point", proved in \cite{CMR}.
In Section~3 we compare various properties of LQEL and CC-manifolds, Proposition\deft\ref{LQELCC}. Examples show that the former class is much larger than the first. In particular, using results by Hwang-Kebekus \cite{HK} (see also \cite{IR}), we see that complete intersections (of high dimension with respect to the multidegree) are CCM, while the only LQELM which are complete intersections are the hyperquadrics, Proposition\deft\ref{3.4}. We also characterize LQEL-manifolds (when $\delta \geq 3$) in terms of the projective geometry of the variety of lines passing through a general point and contained in $X$, Proposition\deft\ref{3.5} and \cite[Theorem 2.3]{QEL}. The main properties of LQEL-manifolds are collected
in Theorem\deft\ref{T}. In particular, their classification is reduced to the case when the Picard group is generated by the class of the hyperplane section. Conjecturally, QEL-manifolds with $\Pic(X) \cong \mathbb{Z} \langle \mathcal {O}_X(1)\rangle$ are linear sections of rational homogeneous manifolds (which are completely classified), see Remark\deft\ref{3.8}. The last section contains two applications. The first one, due to B. Fu \cite{Fu}, gives a substantial improvement of one of the main results from \cite{KS}
 and is based on ideas developped in \cite{QEL}, \cite{IR} and the present paper. The other application is a new proof of the classification, due to L. Ein \cite{Ein}, of manifolds with small dual. To the best of our knowledge, it is for the first time that classification of dual defective manifolds is connected to (and, conjecturally,  even reduced to, cf. Remark\deft\ref{4.5}) that of secant defective ones, via the key concept of (L)QELM.

\section{Preliminaries}

We work over the field of complex numbers. Notation and terminology are the same as in \cite{QEL}; we recall below some of the relevant facts.  

Let $X\subset \Proj^N$ be an irreducible non-degenerate projective
variety of dimension~$n$. Let
\[SX=\overline{\textstyle{\bigcup\limits_{{\scriptstyle {x\neq y}}\atop{\scriptstyle{  x,y\in
X}}}}\langle x,y\rangle }\subseteq\p^N\] be the secant variety to $X$; see also
the construction in \eqref{joindiagram}.

Clearly $\dim(SX)\leq\min\{N,2n+1\}$. If $\dim(SX)<2n+1$, then $X\subset\p^N$ is
said to be {\it secant defective}. The {\it secant defect} of
$X\subset\p^N$ is equal to  $\delta(X):=2n+1-\dim(SX)$.

For $p\in SX\setminus X$, the closure of the locus of couples of
distinct  points of $X$ spanning secant lines passing through $p$ is
called {\it the entry locus of $X$ with respect to $p\in SX$} and it
will be denoted by $\Sigma_p(X)$. The closure of the locus of  secant lines
to $X$ passing through $p$ is a cone over $\Sigma_p(X)$, let us call
it $C_p(X)$.  If $X\subset\p^N$ is smooth, then
$\Sigma_p(X)=C_p(X)\cap X$ as schemes for $p\in SX$ general; see, for
example, \cite[Lemma 4.5]{FR}. Moreover it is easy to see that for
$p\in SX$ general,  $\Sigma_p(X)$ is  equidimensional of dimension
equal to $\delta(X)$. Thus $\dim(C_p(X))=\delta(X)+1$. In general
$\Sigma_p(X)$ may be reducible.
\begin{equation}\begin{split} &\mbox{Let $X \subset \p^N$ be an irreducible non-degenerate projective manifold of}\\
 &\mbox{dimension $n$ and secant defect $\delta \geq 0$.} \end{split}\label{eq1.1}
\end{equation}

In the following definition, we consider varieties having the {\it simplest entry locus}.

\begin{Definition}[(cf. also \cite{KS, QEL,IR})] 
Let $X$ be as in \eqref{eq1.1}.  
\begin{enumerate}
\item[(i)]
$X$  is said to be a {\it quadratic entry locus
manifold of type $\delta\geq 0$}, briefly a {\it QEL-manifold of
type $\delta$}, if for general $p\in SX$ the entry locus
$\Sigma_p(X)$ is a quadric hypersurface of dimension
$\delta$.
\item[(ii)]
$X$ is said to be a {\it local quadratic entry locus
manifold of type $\delta\geq 0$}, briefly an {\it LQEL-manifold of type
$\delta$}, if through two general points
there passes a  quadric hypersurface of dimension $\delta$ contained in $X$.
\item[(iii)]
$X$ is
said to be a {\it conic-con\-nected manifold}, briefly a {\it
CC-manifold}, if through two general points of $X$ there passes an irreducible conic contained in $X$. 
\end{enumerate}
\end{Definition}

Note that, for $\delta =0$, being an LQEL-manifold imposes no restriction on $X$. Clearly, any QEL-manifold is LQEL and any defective LQEL-manifold is CC.

\begin{Lemma} [(cf. Lemma 1.2 in \cite{QEL})] \label{1.4} Let $X$ be an LQEL-manifold with  $\delta>0$ and let $x,y \in X$ be general points. There is a unique quadric hypersurface of dimension $\delta$, say $Q_{x,y}$, passing through $x,y$ and contained in $X$. Moreover, $Q_{x,y}$ is irreducible.  
\end{Lemma}

The
following proposition is easy to prove.

\begin{Proposition}\label{projectQEL}
Let $X\subset\p^N$ be an irreducible non-degenerate smooth
projective variety. 
\begin{enumerate}
\item[(i)] If $X'\subset\p^M$, $M\leq
N-1$, is an isomorphic projection of $X$, then $X'$ is an LQEL-manifold if and only if
$X$ is an LQEL-manifold.

\item[(ii)] If $X$ is an (L)QEL-manifold of type $\delta\geq
1$, then  a general hyperplane section
is an (L)QEL-manifold of type
$\delta-1$.

\item[(iii)] If  $X$ is a QEL-manifold and  $SX=\p^N$,
then $X$ is linearly normal.
\end{enumerate}
\end{Proposition}
\begin{proof} (i) and (ii) are standard and left to the reader, so that we shall prove only part (iii). Suppose $X\subset\p^N$ were the isomorphic projection of $\overline X\subset\p^{N+1}$
from a point $q\in\p^{N+1}\setminus S\overline X$. 
Since $N=\dim(SX)=\dim(S\overline X)$, $S\overline X$ would be a hypersurface of degree at least two. A general $p\in SX$ would  be the projection from $q$ of at least two different points
$p_1,p_2$ belonging to $\langle q,p\rangle \cap S\overline X$. Then the projections of the entry loci $\Sigma_{p_1}(\overline X)$ and $\Sigma_{p_2}(\overline X)$ yield two irreducible components of $\Sigma_p(X)$. This is a contradiction, since $\Sigma_p(X)$
is a smooth, hence irreducible, $\delta$-dimensional quadric hypersurface; see 
 \cite[pp.\ 964--966]{FR}.
\end{proof}

If $x \in X \subset \p^N$ is a smooth point, we denote by ${\bf T}_x X$ the affine Zariski tangent space at $x$ and by $T_x X$ its projective closure in $\p^N$.

The dimension of the image of the projection of
$X\subset\p^N$ from a general tangent space to $X$, called the
{\it tangential projection} of $X\subset\p^N$, is easily computed via
Terracini Lemma (see, for example, \cite[Section 1]{QEL}).
 Let $x\in X\subset\p^N$ be a
general point and let
\begin{equation}\label{tangentdef}
\pi_x:X\map W_x\subseteq\p^{N-n-1},
\end{equation}
 be the projection of $X$ from
$T_xX$. We have $\dim(W_x)=n-\delta$, so that a 
general fiber of $\pi_x$ is of pure 
dimension $\delta$. 

\begin{Lemma}[(cf. Lemma 1.6 in \cite{QEL})]\label{2.1}
Let $X\subset \p^N$ be a smooth irreducible non-degenerate variety, and assume $\delta>0$. The irreducible components of the closure of a general fibre of $\pi_x$ are not linear.  
\end{Lemma}

\section{Tangential projections and the geometry of LQEL-manifolds}

The main result of this section is  the following theorem,
generalizing \cite[Corollary\deft4.2]{CMR}, where  the case of
QEL-manifolds of type $\delta=0$ in $\p^{2n+1}$, i.e.\ of {\it
varieties with one apparent double point}, was considered.

\begin{Theorem} \label{rationality} Let $X\subset \Proj^{N}$
be as in \eqref{eq1.1} and let $x\in X$ be a general point. Then:
\begin{enumerate}
\item[(i)] If $X$ is a QEL-manifold of type $\delta\geq 0$,
 the projection from a general $\delta$-codi\-mensional subspace
of $T_xX$ passing through $x$ is birational onto its image.
 In particular, if $SX=\p^N$, 
 then $X$ is rational.

\item[(ii)] Conversely, if the projection from a general $\delta$-codimensional
subspace of $T_xX$ passing through $x$ is birational, then $X$ is an
LQEL-manifold of type~$\delta$.

\item[(iii)] If $X\subset\p^N$ is an LQEL-manifold of type $\delta>0$, then $X$ is rational.
\end{enumerate}
\end{Theorem}

We recall some
results from \cite{CMR} (see also \cite{CR})~in order to point out a relation
between the entry loci of a variety and a general fiber of the tangential projection. We believe this relation, contained in Proposition\deft\ref{entrydegeneration} below, is interesting~in~itself.

Let $X\subset \p^N$ be an irreducible non-degenerate projective variety and let 
\[
S_X:=\overline{\{(x,y,z)\mid x\neq y,\, z\in \langle x,y\rangle \}}\subset
X\times X\times\p^N,\] be the {\it abstract secant variety}  of
$X\subset\p^N$, which is an irreducible projective variety of
dimension $2n+1$. Let us consider the projections of $S_X$ onto the
factors $X\times X$ and $\p^N$,
\begin{equation}\label{joindiagram}\raisebox{.7cm}{\xymatrix{
          &S_X\ar[dl]_{p_1}\ar[dr]^{p_2}       &             \\
  X\times X& & \p^N.} }            
\end{equation}
With this notation we get
\[p_2(S_X)=\overline{\textstyle\bigcup\limits_{{\scriptstyle x\neq y}\atop{\scriptstyle x,y\in X}}\langle x,y\rangle }=SX\subseteq\p^N.\]

 Let $L=\langle x,y\rangle $ with $x\in X$ and
$y\in X$ general points, i.e.\ $L$ is a general secant line to $X$,
and let $p\in\langle x,y\rangle\subseteq SX\subseteq \Proj^N$ be a general point.
We fix coordinates on $L$ so that the coordinate of $x$ is $0$;  let
$U$ be an open subset of $\mathbb{A}^1_{\C}\subset L$ containing
$0=x$. Let $p_2:S_X\to SX\subseteq\p^N$ be as above and let
\[Z_U=p_2^{-1}(U)\subset S_X.\] 
By shrinking up $U$, we can
suppose that $p_2:Z_U\to U$ is flat over $U\setminus \{0\}$ and that
$\dim(Z_U)_t=\delta(X)$ for every $t\neq 0$. The projection of
$p_1((Z_U)_t)$ onto one of the factors is $\Sigma_t$, the entry locus of $X$ with respect to
$t$ for every $t\neq 0$.

Moreover, by definition, a point $(r,s)\in X\times X$, $r\neq s$,
belongs to $(Z_U)_t$, $t\neq 0$, if and only if $t\in \langle r,s\rangle $, that
is if and only if $(r,s)\in p_2^{-1}(t)$. Thus, if $\psi_t:X\map
\Proj^ {N-1}$ is the projection from $t$ onto a disjoint $\p^{N-1}$,
we can also suppose that $\psi_t$ is a morphism for every $t\neq 0$
and a rational map not defined at $x=0$ for $t=0$. The above
analysis says that the  {\it abstract entry locus} $(Z_U)_t$, $t\neq 0$,
can be considered as the closure in $X\times X$ of the double point
locus scheme of $\psi_t$, minus the diagonal $\Delta_X\subset X\times
X$.

Let $T=\langle T_xX,y\rangle $, so that $T$ is a general $\Proj^{n+1}$ containing
$T_xX$ and a general point $y\in X$. By definition
$\pi_x^{-1}(\pi_x(y))=T\cap X\setminus (T_xX\cap X)$, notation as in \eqref{tangentdef}. Let
\[F_y=\overline{\pi_x^{-1}(\pi_x(y)}),\]
be the closure of the fiber of $\pi_x$ through $y$. Every
irreducible component of $F_y$ has dimension $\delta(X)$ by
Terracini Lemma and by the generality of $y$, see the discussion
after \eqref{tangentdef}. Generic smoothness ensures also that there
exists only one irreducible component of $F_y$ through $y$.

By using the same ideas as in \cite{CMR}, we have the following result,
not explicitly stated in loc.\ cit., because a slightly different
degeneration was considered. For more details about the construction
recalled above and below, we refer to \cite[Sections 3 and 4]{CMR} and 
\cite[Section\deft2]{CR}.

\begin{Proposition}\label{entrydegeneration}
Let notation be as above. The closure of the fiber of $\pi_x$
through a general point $y\in X$ is contained in the flat limit of
the family $\{(Z_U)_t\}_{t\neq 0}$. In other words, the closure of a
general fiber of the tangential projection is a degeneration of the
general entry~locus~of~$X$.
\end{Proposition}
\begin{proof} We shall look at $\psi_t$ as a family of morphisms
and study the limit of the {\it double point scheme}~$(Z_U)_t$.

 Consider the products ${\mathcal X}=X\times U$
and $\Proj_U= \Proj^ {N-1}\times U$. The projections $\psi_t$, for
$t\in U $, fit together to give a rational map $\psi: {\mathcal
X}\map \Proj_U$, which is defined everywhere except at the pair
$(x,0)$. In order to extend the projection not defined at
$x\in X$, we have to blow up ${\mathcal X}$ at $(x,0)$. Let $\sigma:
\widetilde{\mathcal X}\to {\mathcal X}$ be this blowing-up and let
$Z\simeq \Proj ^{n}$ be the exceptional divisor. Looking at the
obvious morphism $\phi: \widetilde{\mathcal X}\to U$, we see that
this is a flat family of varieties over $U$. The fiber
$\widetilde{\mathcal X}_t$ over a point $t\in U\setminus\{0\}$ is
isomorphic to $X$, whereas the fiber $\widetilde{\mathcal X}_0$ over
$t=0$ is of the form $\widetilde{\mathcal X}_0=\widetilde{X}\cup Z$,
where $\widetilde{ X}\to X$ is the blowing-up of $X$ at $x$, and
$\widetilde{X}\cap Z=E$ is the exceptional divisor of this blowing-up,
the intersection being transverse. Reasoning as in \cite[Lemma
3.1]{CMR}, it is easy to see that $\psi_0$ acts on $\widetilde{X}$ as
the projection from the point $0=x$, while it maps $Z$
isomorphically onto the linear space $\psi_0(T)=\p^n$. This
immediately implies that every point of $T\cap X$, different from
$x$, appears in the ``double point scheme" of $\psi_0:\widetilde{X}\cup
Z\to\p^{N-1}$. Therefore $F_y$, being of dimension $\delta(X)$, is
contained in the flat limit of $\{(Z_U)_t\}_{t\neq 0}$, proving the
assertion.
\end{proof}

For an irreducible variety $X\subset\p^N$ we denote by
$\mu(X)$ the number of secant lines passing through a general point
of $SX$. If $\delta(X)>0$, then $\mu(X)$ is infinite, while for
$\delta(X)=0$ the above number is finite and in this case
\[\nu(X)=\mu(X)\cdot\deg(SX)\]
is called {\it the number of apparent double points} of
$X\subset\p^N$. With these definitions we obtain the following
generalization of \cite[Theorem 4.1]{CMR} (see also 
\cite[Theorem\deft2.7]{CR}).

\begin{Theorem}\label{prebronowski} Let $X\subset \Proj^{N}$ be as in \eqref{eq1.1}.
If $\delta(X)=0$, then
\[0<\deg(\pi_x)\leq\mu(X).\]
In particular for a QEL-manifold of type $\delta=0$, the general
tangential projection is birational.

If $X\subset\p^N$ is a QEL-manifold of type $\delta>0$, then the
general fiber of $\pi_x$ is irreducible. More precisely the closure
of the fiber of $\pi_x$ passing through a general point $y\in X$ is
the entry locus of a general point $p\in\langle x,y\rangle $, i.e.\ a smooth
quadric hypersurface.
\end{Theorem}

\begin{proof} If $\delta(X)=0$, then for $t\in U\setminus \{0\}$ the
$0$-dimensional scheme $(Z_U)_t$ has length~equal to $2\mu(X)$. The
$0$-dimensional scheme $F_y$ contains $\deg(\pi_x)$ isolated points,
yielding $2\deg(\pi_x)$ points in the flat limit of
$\{(Z_U)_t\}_{t\neq 0}$ by Proposition \ref{entrydegeneration} and 
proving the first~part.

 Suppose $X$ is a QEL-manifold of type $\delta>0$. Then
for every $t\neq 0$ the $\delta$-dimensional scheme $(Z_U)_t$ is a
smooth quadric hypersurface by definition of QEL-manifolds. The
fiber $F_y$ contains the entry locus $\Sigma_p$ of a
general point $p\in \langle x,y\rangle $, which is a smooth quadric hypersurface
of dimension $\delta$ passing through $x$ and $y$. By proposition
\ref{entrydegeneration} the variety $F_y$ is also contained in the flat
limit of $\{(Z_U)_t\}_{t\neq 0}$.  Therefore $F_y$
coincides with $\Sigma_p$. In fact, in this case the
family $\{(Z_U)_t\}_{t\neq 0}$ is constant.
\end{proof}

\begin{proof}[Proof of Theorem {\rm \ref{rationality}}]  Suppose $X\subset\p^N$ is a QEL-manifold of
type $\delta\geq 0$. If $\delta=0$ then the
first part of Theorem \ref{prebronowski} yields that
$\pi_x$ is birational onto its image (see also \cite[Corollary 4.2]{CMR}).

 Suppose from now on $\delta>0$. The projection from
a general codimension $\delta$ linear subspace $L\subseteq T_xX$
passing through $x$ is a rational map
$\pi_L:X\map\p^{N-n+\delta-1}$. 
For general $y\in X$, the linear space $\langle L,y\rangle$ is obtained by cutting $\langle T_xX,y\rangle$ with $\delta$
general hyperplanes $H^y_1,\ldots, H^y_{\delta}$ passing through $x$ and $y$. From the second part of Theorem \ref{prebronowski}, it follows
that, for general $y\in X$, $\overline{\pi_x^{-1}(\pi_x(y))}$ is a $\delta$-dimensional quadric $Q_{x,y}$. Then $$\pi_L^{-1}(\pi_L(y))\subseteq Q_{x,y}\cap H^y_1\cap\ldots\cap H^y_{\delta}=\{x,y\},$$
yielding the birationality of $\pi_L$ onto its image.
If $SX=\p^N$, then $N-n+\delta -1=n$. 
Part (i) is now completely
proved.

Suppose we are in the hypothesis of (ii). Let
$L=\p^{n-\delta}\subset T_xX$ be a general linear subspace passing
through $x$. $L$ is the tangent space of a general
codimension $\delta$ linear section of $X\subset\p^N$ passing
through $x$, let us say $Z$. Thus the restriction of $\pi_L$ to $Z$, 
$\pi_{L|Z}:Z\map \p^{N-\delta-\dim(Z)-1}=\p^{N-n-1}$, is
the projection from $L=T_xZ$. Since $\pi_L$ restricted to $X$ is
birational onto its image, also $\pi_{L|Z}$ is easily seen to be
birational onto its image. Moreover, looking at $\pi_{L|Z}$ as the
projection from $T_xZ$, we get
$\pi_{L|Z}(Z)=\pi_x(X)=W_x\subseteq\p^{N-n-1}$. Each irreducible component of a general
fiber of $\pi_x$ produces at least one point in a general fiber of $\pi_{L|Z}$.  Hence 
$\pi_x:X\map W_x\subseteq\p^{N-n-1}$ has irreducible general fibers of
dimension $\delta$ by the birationality of $\pi_{L|Z}$.

We shall prove inductively that $F_y=\overline{\pi_x^{-1}(\pi_x(y)})$ is a $\delta$--dimensional quadric for general
$y\in X$.
So, there is no loss of generality in supposing $\delta=1$, by
passing to a general linear section; see Proposition~\ref{projectQEL}.
We claim that set theoretically $L\cap F_y=\{x\}$. We have  $F_y\subset
\langle T_xX,y\rangle $ so that $T_xX\cap F_y$ consists  of a finite number of  points.
By the generality of $L$  we
get $L\cap F_y\subseteq \{x\}$. Let $t=\pi_x(y)$. Since $\langle L,t\rangle $ is a
hyperplane in $\langle T_xX,y\rangle =\langle T_xX,t\rangle $, intersecting $\pi_x^{-1}(\pi_x(y))$ transversally
at a unique point, we get that either we are in the case of the
claim or $F_y$ is a line. This last case is excluded by Lemma\deft\ref{2.1}.

Let $q=\pi_{L|Z}^{-1}(t)=\langle L,t\rangle \cap Z\setminus (L\cap Z)$. By
definition
\begin{equation}\label{Lt} \langle L,t\rangle =\langle L,q\rangle .
\end{equation}

Consider the projection from $t$ onto $T_xX$, let us say
$\psi_t:\langle T_xX,t\rangle \map T_xX$. Let $\widetilde{F}_y=\psi_t(F_y)$. By
definition $x\in\widetilde{F}_y$ because $x\in F_y$. Moreover we
claim that $L\cap \widetilde{F}_y$ is supported at $x$, so that
$\widetilde{F}_y$ is a line through $x$. Indeed, if $z\in L\cap
\widetilde{F}_y$, then there exists $w\in\langle z,t\rangle \cap
F_y\subset\langle L,t\rangle \cap F_y=\langle L,q\rangle \cap F_y$, where the last equality
follows from \eqref{Lt}. Thus either $w=x$ or $w=q$ and in any case
$x=\psi_t(w)=z$. Therefore $F_y\subset\langle \widetilde{F}_y,t\rangle \simeq\p^2$;
moreover the Trisecant Lemma implies that the line $\langle x,y\rangle $ cuts transversally $X$, hence also $F_y$, at $x$ and at
$y$. The line $\langle x,y\rangle $ is contained in the plane
$\langle \widetilde{F}_y,t\rangle $, so that
 $\deg(F_y)=2$ and $F_y$ is  a
smooth conic passing through $x$ and $y$, concluding the proof.

Let us prove part (iii). Fix a general point $x\in X$ and denote by $\mathcal{Q}_x$ (an irreducible component of) the family of $\delta$-dimensional 
quadric hypersurfaces contained in $X$ and passing through $x$. Let $\pi: \mathcal{F}_x \to \mathcal{Q}_x$ be the universal family and let $\phi: \mathcal{F}_x \to X$ be the tautological morphism.

Assume first $\delta =1$. Then $\phi$ is surjective by definition of LQEL-manifold and birational by Lemma\deft\ref{1.4}.
Note that $\pi$ has a section, corresponding to the point $x$. Denote by ${\mathcal E}\subset\mathcal{F}_x$ the image of this section.
 Consider the blowing-up $\sigma: X'\to X$ of $X$ at $x$. Since both $\mathcal Q_x$ and $\pi$
are generically smooth, the birational map $\psi=\sigma^{-1}\circ\phi$ is defined at a general point of $\mathcal E$. Moreover, as $\mathcal E$ is contracted by $\phi$ to the point $x$, $\psi$ sends $\mathcal E$ to $E$, the exceptional divisor
of $\sigma$. So we have the following diagram
\begin{equation*}\label{joindiagram1}\raisebox{.7cm}{\xymatrix{
&\mathcal{F}_x  \ar[d]_\pi \ar[dr]^\phi\ar@{-->}[r]^\psi&X'\ar[d]^\sigma\\
&\mathcal{Q}_x&X.
}}
\end{equation*}

 By Zariski Main Theorem the map $\psi^{-1}$
is defined at a  general point of $E$, hence  $\psi$  induces a birational map $\psi_0:\mathcal E\map E=\p^{n-1}$. Thus $\mathcal F_x$
is birational to $\p^n$, being birational to a family of conics with a section over a rational base,  and  $X$ is a rational variety, as claimed. A general version of the above argument appears in  \cite[Proposition\deft3.1]{IN}.

Suppose now that $\delta \geq 2$  and fix $H$ a general hyperplane section of $X$ through $x$. Using Lemma\deft\ref{1.4}, we see that sending a quadric hypersurface through $x$ to its trace on $H$ yields a birational map between the families $\mathcal{Q}_x(X)$ and   $\mathcal{Q}_x(H)$. So, we see inductively that $\mathcal{Q}_x(X)$ is a rational variety of dimension $n-\delta$. Therefore, $\mathcal{F}_x$ is rational, as the family  $\mathcal{F}_x \to \mathcal{Q}_x$ has a section. Being birational to $\mathcal{F}_x$ by Lemma\deft\ref{1.4}, $X$ is rational too.   
\end{proof}

\begin{Remark}\label{Bronowski}  J.\ Bronowski
claims in \cite{Br} that  $X\subset \p^{2n+1}$ is a  variety with
one apparent double point if and only if the projection of $X$ to
$\Proj^n$ from a general tangent space  to $X$ is a birational map;
he also formulates  a generalization to arbitrary secant
$k$-planes, see \cite{CR} for more details. Unfortunately
Bronowski's argument is unclear and we do not know of any convincing
proof  for this statement. \cite[Corollary 4.2]{CMR} (see also
Theorem \ref{rationality}) proves one implication showing the
rationality of varieties with one apparent double point. The open
implication would be a very useful tool for constructing examples.

We can generalize {\it Bronowski Conjecture} to the following: a smooth 
irreducible $n$-di\-men\-sional variety $X\!\subset\!\p^{2n+1-\delta(X)}$  is a
QEL-manifold if and only if the projection~from~a general
codimension $\delta(X)$ linear subspace of $T_xX$ passing through
$x$ is birational.~Theorem\deft\ref{rationality} proves one implication, 
 yielding  the~rationality of QEL-manifolds and
extending \cite[Corollary 4.2]{CMR}. 
One may consult
\cite{CR} for  other generalizations  of
the above conjecture to higher secant varieties.

It is worth mentioning  that the above results reveal the following
interesting picture for the tangential projections of
QEL-manifolds of type $\delta\geq 0$ with $SX=\p^N$: for
$\delta=0$ we project from the whole space and we have varieties
with one apparent double point; at the other extreme we found the
stereographic projection of quadric hypersurfaces, the only
QEL-manifolds of type equal to their dimension.
\end{Remark}

\section{LQEL versus CC}

We recall the following definition from \cite{BaBeIo}.

\begin{Definition}[(\cite{BaBeIo})]\label{bbi} A smooth rational curve $C \subset X$, where $X$ is a projective manifold of dimension $n$, is a {\it quasi-line} if $N_{C|X}\simeq \bigoplus\limits^{n-1}_1 \mathcal{O}_{\p^1}(1).$
\end{Definition}

The relation between the notions  LQEL-manifold and 
CC-manifold is clarified in  the following.

\begin{Proposition}\label{LQELCC}
 Let $X\subset\p^N$ be a CC-manifold of secant defect $\delta$.  Let $C=C_{x,y}$ be
 a general conic through the general points $x,y\in X$ and let $c$ be the point representing $C$ in the Hilbert scheme  of $X$. Let $\mathcal{C}_x$ be the unique irreducible component of the Hilbert scheme of conics passing through $x$ which contains the point $c$. 
\begin{enumerate}
\item[(i)] We have
 $n+\delta\geq -K_X\cdot C= \dim (\mathcal{C}_x)+ 2 \geq n+1$.

\item[(ii)] The equality  $-K_X\cdot C=n+\delta$  holds if and only if $X\subset\p^N$ is
an LQEL-manifold.

\item[(iii)] The following conditions are equivalent:
\begin{enumerate}
\item[(a)] $\dim (\mathcal{C}_x) = n-1$;
\item[(b)] $ C$ is a quasi-line;
\item[(c)] all conics through $x,y$ are non-degenerate. 
\end{enumerate}

\item[(iv)] If $\delta\geq 3$, then $X\subset\p^N$ is a Fano
manifold with $\Pic(X)\simeq\mathbb{Z}\langle \O_X(1)\rangle$ and  index $i(X)=\frac{\dim(\mathcal{C}_x )}2 +1$.
\end{enumerate}
\end{Proposition}
\begin{proof}
We have the universal family $g:\mathcal{F} _x\to \mathcal{C}_x$ and the tautological morphism $f: \mathcal{F}_x \to X$, which is surjective. 
 Since $C\in \mathcal{C}_x$ is a general conic and since
$x\in X$ is general, we get $\dim(\mathcal{C}_x)=-K_X\cdot C-2$.

Take a general point $y\in X$ and a general $p\in \langle x,y\rangle $. The conics
passing through $x$ and $y$ are parameterized by $g(f^{-1}(y))$,
which has pure dimension
\[\dim({\mathcal{F}}_x)-n=\dim(\mathcal{C}_x)+1-n=-K_X\cdot C-1-n.\]
We claim that the locus of conics through $x$ and $y$, denoted by $\mathcal{L}_{x,y}$, has dimension
$\mbox{$-K_X\cdot C-n$}$ and is clearly contained in the irreducible
component of the entry locus (with respect to $p$) through $x$ and $y$. 
Indeed, conics through $x,y$ and another general point $z\in \mathcal{L}_{x,y}$ have to be finitely many.
Otherwise, their locus would fill up the plane $\langle x,y,z\rangle$ and this would imply that the line $\langle x,y\rangle$ is contained in $X$. But we have excluded linear spaces from the definition of CC and LQEL-manifolds. 
Therefore
$\delta\geq -K_X\cdot C-n$, that is $-K_X\cdot C\leq
n+\delta$.

 The locus of conics through $x$ and $y$ is  contained in
$\langle T_xX,y\rangle \cap\langle x,T_yX\rangle $, which is a linear space of dimension
$\delta+1$. Indeed, by Terracini Lemma, 
\[\dim (T_xX\cap T_yX)=2n-\dim (SX)=\delta-1.\] By the Trisecant Lemma $x\not\in T_yX$
and $y\not\in T_xX$, so that \[\dim( \langle T_xX,y\rangle \cap\langle x,T_yX\rangle)=\delta+1.\]

 If $-K_X\cdot C= n+\delta$, then, for $p\in \langle x,y\rangle $ general, the
irreducible component $\Sigma^p_{x,y}$ of the entry locus passing
through $x$ and $y$ coincides with the locus of conics through $x$
and $y$. Thus $\Sigma^p_{x,y}$ is a quadric
hypersurface by the Trisecant Lemma and by the generality of $x$ and
$y$ (if $\delta=n$, $X\subset\p^{n+1}$ is a
quadric hypersurface). So, (ii) is proved.

Next we see (iii). (a) and (b) are equivalent, since the normal bundle of $C$ in $X$ is ample, of degree
$\dim (\mathcal{C}_x)$. 

Assume that (a) holds. The dimension of the subfamily consisting of reducible conics from $\mathcal{C}_x$  has dimension at most $n-2$. Hence their locus is of dimension at most $n-1$ and does not contain the general point $y$. So we have (c).  

Assume that (c) holds. Then, by bend and break there are finitely many conics through $x$ and $y$, giving (a).

 Finally, (iv) follows from  the Barth--Larsen Theorem, the fact that $X$ contains moving conics and (i).
\end{proof}

\begin{Examples}\label{1.9}  (i) 
For $n\geq 3$, let $X$ be a smooth cubic hypersurface in $\p^{n+1}$ or the smooth complete intersection of two hyperquadrics in $\p^{n+2}$. Use e.g.\ \cite[Theorem\deft3.2]{BaBeIo} and induction on $n$ to see that $X$ is conic-connected. In the first case, $X$ is a hypersurface of degree $3$ so it cannot be an LQEL-manifold. In the second case, $\delta(X)=n-1$; as the Picard group is generated by the hyperplane section, $X$ cannot contain hyperquadrics of dimension $n-1$.  
\vskip4pt

(ii) CC-manifolds $X\subset\p^{N+1}$ of secant defect $\delta(X)=\delta-1\geq 2$,
constructed from QEL-manifolds $Z\subset\p^N$ of type $\delta\geq
3$.
\vskip4pt

Let $Z\subset\p^N$ be a QEL-manifold of type $\delta\geq 3$.
Consider $\p^N$ as a hyperplane in $\p^{N+1}$,~take
$q\in\p^{N+1}\setminus\p^N$ and let $W=C_q(Z)\subset\p^{N+1}$ be the
cone over $Z$ of vertex $q$. Let $X=C_q(Z)\cap Q\subset\p^{N+1}$,
where $Q\subset\p^{N+1}$ is a general quadric hypersurface. Then
$SX=C_q(SZ)$, yielding $\delta(X)=\delta(Z)-1\geq 2$. Moreover,
$\Pic(X)\simeq\mathbb{Z}\langle \O_X(1)\rangle $ and $X\subset\p^{N+1}$ is a
CC-manifold of index $i(X)=\frac{n+\delta(X)-1}{2}$.

Indeed, let $\pi:T=\p(\O_Z\oplus\O_Z(1))\to Z$ and let $E\subset T$
be the section at infinity of $\pi$. If
$\phi:T=\p(\O_Z\oplus\O_Z(1))\to\p^{N+1}$ is the tautological
morphism given by (a sublinear system of) $|\O_T(1)|$, then
$\phi(T)=W$ and $X\subset\p^{N+1}$ can be naturally thought of  as an
element of $|\O_T(2)|$. Recall that  $\phi$ restricts to an
isomorphism between $T\setminus E$ and $W\setminus \{q\}$. By adjunction
we get that $\omega_X=\O_X(1-i(Z))$, that is
\[i(X)=i(Z)-1=\frac{n+\delta}{2}-1=\frac{n+\delta(X)-1}{2}.\] In
particular $X\subset\p^{N+1}$ is a Fano manifold. Adjunction formula
also says that the double cover $\pi_q:X\to Z$, induced by the
projection from $q$, $\pi_q:\p^{N+1}\setminus \{q\}\map\p^N$, is
ramified along a hyperplane section of $Z\subset\p^N$.

Take two general points $x,y\in X$ and let $x'=\pi_q(x)$ and
$y'=\pi_q(y)$. Through $x'$ and $y'$ there passes a smooth quadric
hypersurface $Q'\subset\p^{\delta+1}$ of dimension $\delta\geq 3$.
Let $\widetilde{Q}=C_q(Q')\subset\p^{\delta+2}$ and let
$\widetilde{X}=C_q(Q')\cap Q\subset\p^{\delta+2}$. Then
$\widetilde{X}\subset\p^{\delta+2}$ is a (smooth) complete
intersection of two quadric hypersurfaces passing through $x$ and
$y$, so that $X\subset\p^{N+1}$ is a CC-manifold since
$\widetilde{X}\subset\p^{\delta+2}$ is conic-connected by (i).
\vskip4pt

(iii) When $\delta =1$, CC-manifolds and LQEL-manifolds are the same. In this case a general conic through two general points is a quasi-line. Examples are easy to construct, e.g.\ $X=G\cap H_1\cap H_2\cap H_3$, where $G$ is the Grassmannian of lines in $\p^m$, $m\geq 4$, and $H_i$ are general hyperplane sections of its Pl\" ucker embedding.   
\vskip4pt

(iv) Assume $\delta =2$ and $X$ is conic-connected. Let $C$ be a general conic passing through two general points. Then $X$ is an LQEL-manifold if and only if $C$ is not a quasi-line; see 
Proposition~\ref{LQELCC}. Examples of the LQEL-case are given by taking $X=G\cap H_1\cap H_2$, where $G$ and $H_i$ are as above. Examples of the case where $X$ contains quasi-lines are got by applying the construction in (ii) above starting with $Z=G\cap H$, a hyperplane section of the same Grassmannian. 
\end{Examples}
Many examples of CC-manifolds which are not LQEL come from the following:

\begin{Proposition}[(cf.\  Corollary~2.5 in \cite{IR})] $\mbox{    }$\label{3.4}

\begin{enumerate}
\item[(i)] 
If $X\subset \p ^{n+r}$ is a smooth non-degenerate complete intersection of multi-degree $(d_1, d_2, \ldots, d_r)$ 
with $n> 3\big( \sum_1^r d_i -r-1\big)$, then $X$ is a CC-manifold.  

\item[(ii)] If $X$ is a secant defective LQEL-manifold and a complete intersection, then $X$ is a hyperquadric.
\end{enumerate}
\end{Proposition}
\begin{proof} (i) is exactly \cite[Corollary~2.5]{IR}. 

We show (ii). If $n=2$ use \cite[Proposition~3.3]{QEL} to conclude. If $n\geq 3$, by Lefschetz Theorem, we have $\Pic(X) \cong \mathbb{Z}\langle \mathcal{O} _X(1)\rangle$ and Proposition~\ref{LQELCC} gives $i(X)=\frac {n+\delta}2$. Let $r=N-n$. We have 
$2n+1 -\delta = \dim (SX) \leq N = n+r$, so $\delta \geq n+1-r$. Assuming $X$ to be a complete intersection of type $(d_1, \ldots, d_r)$, with $d_i\geq 2$ for all $i$, we get
$i(X)=n+r+1-\sum^r_1 d_i $, hence 
\[ n+2 r +2 -2 \sum^r_1 d_i = 2 i(X) -n = \delta \geq n+1 -r \quad \mbox{and} \quad 4r \leq 2\sum^r_1 d_i \leq 3r+1,\]
so $r=1$.
\end{proof}

The proof of the following criterion for recognizing LQEL-manifolds illustrates the role of CC-manifolds. Recall from \cite{QEL} that, if $x\in X$ is a general point, we denote by $Y_x\subset \p (\mathbf{T}^*_x X)$
the variety  of lines through $x$, contained in $X$.  
 
\begin{Proposition}\label{3.5}
Let $X\subset \p^N$ be as in \eqref{eq1.1} and  assume $\delta \geq 3$. The following assertions are equivalent:
\begin{enumerate}

\item[(i)] $X$ is an LQEL-manifold;

\item[(ii)] if $x\in X$ is a general point, $\dim (Y _x )\geq \frac{n+\delta}{2} -2$ and  $SY_x = \p^{n-1}$.
\end{enumerate}
\end{Proposition}
\begin{proof} (i) implies (ii) was shown in \cite[Theorem~2.3]{QEL} and, in fact, equality holds in (ii).

Assume (ii). $\delta \geq 3$ implies via Barth--Larsen Theorem that $\Pic(X)\cong \mathbb{Z} \langle \mathcal{O} _X(1)\rangle$. Moreover 
\[   i(X) = \dim (Y_x) +2 \geq \frac{n+\delta} 2 \geq \frac {n+3} 2.\]
So we may apply \cite[Theorem~3.14]{HK} to deduce that $X$ is a CC-manifold. In the notation of Proposition~\ref{LQELCC}, we have $2i(X) = - K_X \cdot C \geq n+\delta $. Therefore, $X$ is an LQEL-manifold by combining (i) and (ii) of Proposition~\ref{LQELCC}. 
\end{proof}

The following proposition is an application of \cite[Theorem~2.3]{QEL}. 

\begin{Proposition}\label{3.6} Assume $X$ is an LQEL-manifold of type $\delta <n$.  

\begin{enumerate}
\item[(i)] If $ \delta > \frac n2$ then $ \delta \leq 6$ and $n \leq 10$. These cases are classified in \cite[Corollary~3.1]{QEL}.

\item[(ii)] If $\delta > \frac n3 $ then $\delta \leq 10$ and $n\leq 26$. These cases may also be classified by the same method. 

\item[(iii)] If $k\geq 4$, $\delta >\frac nk$, then $\delta \leq 2(k+2)$. 
\end{enumerate}
\end{Proposition}
\begin{proof}
See \cite[Corollary 3.1]{QEL} for (i). 

Assume $\frac n3<\delta \leq \frac n2$. Using \cite[Theorem~2.8]{QEL}, we find $\delta \leq 10$ and the following possibilities for the pairs $(n,\delta)$: $(2,1), (4,2), (5,2), (7,3), (8,4), (10,4),$ $(13,5), (14,6), (15,7), (16,8), (25,9)$ and $(26,10)$.

To prove (iii) we proceed by induction on $k\geq 3$. We may suppose that $\delta \geq 10$. Assume we have $\frac nk \geq \delta > \frac n{k+1}$ and $k\geq 4$. In the notation of \cite[Theorem 2.4]{QEL}, we let $X_1=Y_x$, $\delta _1 = \delta (X_1)=\delta -2$, $n_1=\dim(X_1)= \frac {n+\delta} 2-2$. If $\delta_1> \frac{n_1} k $ we get $ \delta \leq 2(k+3)$, completing the induction. If $\delta_1 \leq \frac{n_1} k$ it follows
\[\delta \leq \frac {n+4(k-1)} {2k-1}, \qquad \mbox{so}\quad \frac n {k+1} < \frac {n+4(k-1)}{2k-1}.\]
This gives \[10k \leq \delta k\leq n<\frac {4(k^2 -1)}{k-2},\] so $3k^2 -10k +2 <0$. Therefore $k=3$; a contradiction. 
\end{proof}

Consider the following list of examples of QEL-manifolds.
\[\mbox{\hskip6pt\begin{minipage}{11.5cm}
\begin{enumerate}
\item[(i)] $\nu_2(\p^n)\subset\p^{\frac{n(n+3)}{2}}$.

\item[(ii)] 
The projection of $\nu_2(\p^n)$ from the linear
space $\langle \nu_2(\p^s)\rangle $, where $\p^s\subset\p^n$ is a linear subspace; equivalently
$X\simeq \Bl_{\p^s}(\p^n)$ embedded in $\p^N$ by the linear system
of quadric hypersurfaces of $\p^n$ passing through $\p^s$;
alternatively $X\simeq\p_{\p^r}(\mathcal{E})$ with
$\mathcal{E}\simeq\O_{\p^r}(1)^{\oplus n-r}\oplus\O_{\p^r}(2)$,
$r=1,2, \ldots, n-1$, embedded by $|\O_{\p(\mathcal{E})}(1)|$. Here $N=\frac{n(n+3)}{2}-\binom{s+2}{2}$ and  $s$ is an integer such that $0\leq s\leq n-2$.
\item[(iii)]  A hyperplane section of the Segre embedding 
$\p^a\times\p^{b}\subset\p^{N+1}$. Here $n\geq 3$ and  $N=ab+a+b-1$, where $a\geq 2$ and $b\geq 2$ are    such that $a+b=n+1$.

\item[(iv)] $\p^a\times\p^b\subset\p^{ab+a+b}$ Segre embedded, where $a,b$ are positive integers such that
$a+b=n$.
\end{enumerate}\end{minipage}}\tag{$*$}\]

The essential properties of LQEL-manifolds are collected in the next theorem. It follows by putting together Theorem~\ref{rationality}\,(iii), \cite[Theorem~2.2]{IR}, and  \cite[Theorem~2.3\,(4d)]{QEL}.

\begin{Theorem}\label{T} Let $ X\subset \p^N$ be a defective LQEL-manifold. Then $X$ is Fano and rational. Moreover, either $X$ is an isomorphic projection of one of the manifolds listed in {\rm ($*$)} or $\Pic(X) \cong \mathbb{Z} \langle\mathcal {O}_X(1)\rangle$ and $i(X)= \frac{n+\delta} 2$. If $\delta \geq 3$, $Y_x \subset \p^{n-1}$ is a QEL-manifold of type $\delta -2$, dimension $\frac{n+\delta}2-2$ and such that $SY_x=\p^{n-1}$.  
\end{Theorem}

\begin{Remark} \label{3.8}
Via the above theorem, the classification of secant defective QEL-manifolds is reduced to the case where the Picard group is $\Z$. Let us say that $X\subset  \p^N$ is {\it maximal} if $X$ is not a hyperplane section of some (non-degenerate) manifold $X' \subset \p^{N+1}$. The following tempting conjecture would, if true, lead to a complete classification of defective QEL-manifolds: 
\vskip4pt

{\it Any maximal defective QEL-manifold with Picard group $\Z$ is homogeneous. } 
\vskip4pt

Note that homogeneous manifolds arising from irreducible representations of (semi)simple complex algebraic groups are QEL-manifolds, see \cite[Chapter III]{Zak} and \cite{Kaji}, and the secant defective ones are classified completely, see 
 loc. cit. 
Moreover, the results in \cite[Section\deft3]{QEL} confirm the conjecture for $\delta>\frac n2$. 
The next two finiteness results for QEL-manifolds with $\delta < n$ would follow from the above conjecture:
\begin{enumerate}
\item[(i)] $\delta \leq 8$;
\item[(ii)] if $\delta \geq 5$, then $n\leq 16$.
\end{enumerate}
Proposition~\ref{3.6} may be seen as supporting these finiteness expectations. Also, in \cite[Corollary 2]{Fu}, it is proved that $\delta \leq \frac{n+8}3$ holds for any LQEL-manifold with $\delta <n$. In particular, (i) follows from (ii).
\end{Remark}

\section{Two applications}

The first application is due to B. Fu \cite{Fu}, who found it relying on the ideas and techniques from \cite{QEL}, \cite{IR} and Proposition\deft\ref{LQELCC}. We mention his result in order to illustrate the usefulness of our point of view. 
\begin{Theorem}[(cf.\ Theorem~2 in \cite{Fu})] Let $X$ be as in \eqref{eq1.1} and assume that $X$ is swept out by hyperquadrics of dimension greater than 
$\big[\frac n2\big]+1$, all passing through a fixed point $x\in X$.  Then $X$ is a hyperquadric. 
\end{Theorem}

This result substantially improves the main application in \cite{KS}, with a much shorter proof. 
\vskip6pt

Our second application concerns the classification of manifolds with small duals. 
 
For an irreducible variety $Z\subset\p^N$, we define
$\defect(Z)=N-1-\dim(Z^*)$ as the  dual defect of $Z\subset\p^N$,
where $Z^*\subset\p^{N*}$ is the dual variety of $Z\subset\p^N$. In
\cite[Theorem 2.4]{Ein} it is proved that if $\defect(X)>0$, then
$\defect(X)\equiv n (\mod 2)$, a result usually attributed to
Landman.  Moreover, Zak Theorem on Tangencies implies that
$\dim(X^*)\geq\dim(X)$ for a smooth non-degenerate variety
$X\subset\p^N$; see \cite[I.2.5]{Zak}.

  We combine the geometry  of CC
and LQEL-manifolds to give a  new proof
of  \cite[Theorem 4.5]{Ein}.
 Our approach 
avoids the use of Beilinson spectral sequences
and more sophisticated computations as in \cite[4.2,~4.3,~4.4]{Ein}.

We begin by recalling some basic facts from \cite{Ein}.

\begin{Proposition}\label{scrolls}
Let $X\subset\p^N$ be as in \eqref{eq1.1} and assume that $\defect(X)>0$. Then
\begin{enumerate}
\item[(i)] {\rm (\cite[Theorem 2.4]{Ein})}  Through a general point $x\in X$
there passes a line $L_x\subset X$ such that $-K_X\cdot
L_x=\frac{n+\defect(X)+2}{2}$, so that $\defect(X)\equiv n$ (\mod 2);
\item[(ii)]{\rm (\cite[Theorem 3.2]{Ein})} $\defect(X)=n-2$ if and only if $X\subset\p^N$ is a scroll
over a smooth curve, i.e.\ it is a $\p^{n-1}$-bundle over a smooth
curve, whose fibers are linearly embedded.
\end{enumerate}
\end{Proposition}

The following proposition reinterprets the result of Hwang and
Kebekus \cite[Theorem~3.14]{HK} on Fano manifolds with large index.

\begin{Proposition}[(cf.\ \cite{HK}, see also Proposition 2.4 in \cite{IR})]\label{dualLQEL} Let $X\subset\p^N$ be as in \eqref{eq1.1}. Assume that $X$ is a Fano manifold with
$\Pic(X)\cong\mathbb{Z}\langle \O_X(1) \rangle $ and let $x\in X$ be a general point.
\begin{enumerate}

\item[(i)] If
$i(X)>\frac{2n}{3}$, then $X$ is a CC-manifold with
$\delta> \frac{n}{3}$.

\item[(ii)] If $\defect(X)>0$ and  
 $\defect(X)>\frac{n-6}{3}$,  then $X$ is a
CC-manifold with $\delta\geq
\defect(X)+2$. Moreover, if $\delta=\defect(X)+2$, then  $X$ is an
LQEL-manifold of type $\delta=\defect(X)+2$.
\item[(iii)] If $X$ is  an LQEL-manifold of type $\delta$
and $\defect(X)>0$, then $\delta=\defect(X)+2$.
\end{enumerate}
\end{Proposition}
\begin{proof} 
See \cite{IR} loc.cit.\ for the proof of (i).

In the hypothesis of (ii), \cite[Theorem 2.4]{Ein}  yields
$i(X)=\frac{n+\defect(X)+2}{2}>\frac{2n}{3}$ so that
$X$ is a CC-manifold by (i).
 Proposition~\ref{LQELCC} yields $\delta\geq\defect(X)+2$ and also the remaining assertions
of (ii) and (iii).
\end{proof}

We recall that according to Hartshorne Conjecture, if
$n>\frac{2}{3}N$, then $X\subset\p^N$ should be a complete
intersection and that complete intersections have no dual defect.
Thus, assuming Hartshorne Conjecture, the following result yields the complete list of manifolds
$X\subset\p^N$ such that $\dim(X^*)=\dim(X)$. The second
part says that under the LQEL hypothesis the same results hold without any
restriction (see also Remark~\ref{4.5} below).

\begin{Theorem}\label{Eindual} Let $X\subset\p^N$ be as in \eqref{eq1.1} and assume that $\dim(X)=\dim(X^*)$.
\begin{enumerate}
\item[(i)] {\rm (\cite[Theorem 4.5]{Ein})} If $N\geq\frac{3n}{2}$, then $X$ is projectively
equivalent to one of the following:
\begin{enumerate}
\item[(a)] a smooth  hypersurface $X\subset\p^{n+1}$, $n=1, 2$;

\item[(b)] a Segre variety $\p^1\times\p^{n-1}\subset\p^{2n-1}$;

\item[(c)] the Pl\" ucker embedding $\mathbb{G}(1,4)\subset\p^9$;

\item[(d)] the $10$-dimensional spinor variety
$S^{10}\subset\p^{15}$.
\end{enumerate}
\item[(ii)] If $X$ is an LQEL-manifold, then it is
projectively equivalent either to  a smooth quadric hypersurface
$Q\subset\p^{n+1}$ or to a variety as in {\rm (b), (c), (d)} above.
\end{enumerate}
\end{Theorem}
\begin{proof} Clearly $\defect(X)=0$ if and only if
$X\subset\p^{n+1}$ is a hypersurface, giving case (a), respectively that of 
quadric hypersurfaces. From now on we suppose
$\defect(X)>0$ and hence  $n\geq 3$. By parts (i) and (ii) of Proposition~\ref{scrolls}, $\defect(X)=n-2$ and $N=2n-1$ if and only if we are in
case  (b); see also \cite[Theorem 3.3, c)]{Ein}.

Thus, we may assume $0<\defect(X)\leq n-4$, that is
$N\leq 2n-3$. Therefore  $\delta\geq 4$ and
  $X$ is a Fano manifold with
$\Pic(X)\cong\mathbb{Z}\langle \O_X(1) \rangle $.  
Moreover, in 
case (i), $\defect(X)=N-n-1>\frac{n-6}{3}$ by hypothesis. Thus Proposition~\ref{dualLQEL} yields that $X$ is also a  CC-manifold with 
$\delta\geq\defect(X)+2$.
Taking into account
also the last part of  Proposition~\ref{dualLQEL}, from now on we can
suppose that $X$ is a CC-manifold with
$\delta\geq\defect(X)+2\geq 3$.

We have
$n-\delta\leq N-1-n =\defect(X)\leq\delta-2$, that is
$\delta\geq\frac{n}{2}+1$. Zak Linear Normality Theorem implies
$SX=\p^N$, so that \[N=\dim(SX)=2n+1-\delta\leq \frac{3n}{2}.\]

Since $N\geq\frac{3n}{2}$, we get  $N=\frac{3n}{2}$,
$\delta=\frac{n}{2}+1=\defect(X)+2$ and $n$  even. Therefore
$X$ is an LQEL-manifold of type $\delta=\frac{n}{2}+1$ by
Proposition~\ref{dualLQEL}. \cite[Corollary~3.1]{QEL}  concludes the
proof, yielding cases (c) and (d).
\end{proof}
\begin{Remark}\label{4.5} Let $X\subset \p^N$ be as in \eqref{eq1.1}. Assume that $\defect(X) >0$ and $\Pic (X)\cong \mathbb{Z}\langle \mathcal{O}_X(1)\rangle$. We conjecture that $X$ is an LQEL-manifold and even a QEL-manifold if moreover assumed to be linearly normal. Combined with Remark~\ref{3.8}, this would imply that maximal dual defective manifolds with $\Pic (X)\cong \mathbb{Z}\langle \mathcal{O}_X(1)\rangle$ are homogeneous, as already conjectured in \cite{BS}. 
\end{Remark}

Our last application is the following:
 
\begin{Proposition} Let $X\subset \p^N$ be as in \eqref{eq1.1}, with $\Pic(X) \cong \mathbb{Z} \langle \mathcal{O}_X(1)\rangle$ and $\defect(X)> 0$. Assume moreover that  $\defect (X) >\frac {n-6} 3 $ and $N\geq \frac {5n+2}3$. Then $SX\neq \p^N$.
\end{Proposition}

\begin{proof}
By Proposition \ref{dualLQEL}, $X$ is a CC-manifold. Assume that $SX=\p^N$. We get
\[N=\dim (SX) = 2n+1-\delta \geq \frac {5n+2}3,\]
so $\delta \leq \frac{n+1}3$. It follows that $3 \defect(X) \geq n-5 \geq 3\delta -6$. By Proposition~\ref{dualLQEL}, $X$ is an LQEL-manifold of type $\delta= \frac{n+1}3$ and $\defect(X)=\delta-2$. From the list in the proof of Proposition~{\rm \ref{3.6}(ii)} it follows that $\delta \leq 2$, so $\defect(X)=0$. This is a contradiction.

\end{proof}

{\bf Acknowledgements}.
Both authors are grateful to the organizers of the Conference ``Projective varieties with unexpected properties", that took place in Siena, between 8--13 of June 2004. Our collaboration started while taking part in this very pleasant and fruitful mathematical event. 

\end{document}